\documentclass[twoside]{article}
\usepackage{blindtext} % Package to generate dummy text throughout this template

\usepackage[sc]{mathpazo} % Use the Palatino font
\usepackage[T1]{fontenc} % Use 8-bit encoding that has 256 glyphs
\usepackage{microtype} % Slightly tweak font spacing for aesthetics

\usepackage[english]{babel} % Language hyphenation and typographical rules

\usepackage[hmarginratio=1:1,top=32mm,columnsep=20pt]{geometry} % Document margins
\usepackage[hang, small,labelfont=bf,up,textfont=it,up]{caption} % Custom captions under/above floats in tables or figures
\usepackage{booktabs} % Horizontal rules in tables
\usepackage{graphicx}
\usepackage{lettrine} % The lettrine is the first enlarged letter at the beginning of the text

\usepackage{enumitem} % Customized lists
\setlist[itemize]{noitemsep} % Make itemize lists more compact

\usepackage{abstract} % Allows abstract customization
\usepackage{titlesec} % Allows customization of titles
\renewcommand\thesection{\Roman{section}} % Roman numerals for the sections
\renewcommand\thesubsection{\roman{subsection}} % roman numerals for subsections
\titleformat{\section}[block]{\large\scshape\centering}{\thesection.}{1em}{} % Change the look of the section titles
\titleformat{\subsection}[block]{\large}{\thesubsection.}{1em}{} % Change the look of the section titles

\usepackage{fancyhdr} % Headers and footers
\usepackage{amsmath}
\usepackage{titling} % Customizing the title section

\usepackage{hyperref} % For hyperlinks in the PDF

\pretitle{\begin{center}\Huge\bfseries} % Article title formatting
\posttitle{\end{center}} % Article title closing formatting

\title{Problem of optimal control for bilinear systems with endpoint constraint}
\author{%
\textsc{ }\\
\normalsize S. Yahyaoui$^1$,   L. Lafhim$^2$ and M. Ouzahra$^1$\\
(1) M2PA Laboratory,  Department
of Mathematic and Informatics \\
ENS. University of Sidi Mohamed Ben Abdellah,  Fez, Morocco. \\
(2) LASMA Laboratory, FSDM, University of Sidi Mohamed Ben\\ Abdellah,  Fez, Morocco.
\normalsize \href{} \normalsize \href{}{}% Second author's email address
}
\date{}
\renewcommand{\maketitlehookd}{%
\begin{abstract}
In this work, we will investigate the question of optimal control for bilinear systems with constrained endpoint. The optimal control will be characterized through a set of unconstrained minimization problems that approximate the former. Then a class of bilinear systems for which the optimal control can be expressed as a time-varying feedback law  will be identified. Finally, applications to parabolic and hyperbolic partial differential equations are provided.\\

       \textbf{ Key words:} Quadratic cost, optimal control, endpoint constraint, bilinear systems.
\end{abstract}
}

\newtheorem{thm}{Theorem}
\newtheorem{lem}[thm]{Lemma}
\newtheorem{rem}[thm]{Remark}

\newtheorem{df}[thm]{Definition}

\begin{document}
\maketitle
\section{Introduction and the problem statement}
Linear systems are usually preferable when approximating nonlinear dynamical processes for their simplicity. However, there are many other practical situations for which  bilinear models are more appropriate (see \cite{bea11,Bradly,Alami,Khapalov,Mohler,Wei} and the references therein).  In general, a problem of control aims to achieve a certain degree of performance for the system at hand using suitable control laws among available options. If this is indeed feasible, then one usually aims to achieve this performance while optimizing a certain criterion. A problem of optimal control is an optimization problem on a reasonable set described by dynamic constraints.  As an interesting example, the question of describing the best control among those that allow to reach a  desired state with minimal cost or energy.	 Such problems arise in various applications, such as the optimization of hydrothermal systems and non-smooth modeling in mechanics and engineering, etc.  (see e.g. \cite{bayon,bayon1,cesari,demyanov,lopes}).
		The problem of  optimal control for  bilinear and semi-linear systems with unconstrained endpoint has been treated by many authors (see \cite{Bradly,cannarsa 92,Alami,Li and Young,linag,Boukhari R,Boukhari 2017}). The  question  of  optimal control with endpoint constraint has been treated  in the context of linear and semi-linear systems with additive controls   (see \cite{Fattorini,Li and Young} and the references therein).
The  approach is based on the Pontryagin's maximum principle. The main goal of this paper is to study a quadratic optimization problem with a restricted endpoint state.	
In the case of a bounded set of admissible control, we will characterize the optimal control either  for exactly or approximately attainable states. This problem can be  formulated  as an optimization problem with endpoint constraint, which can also be approximated by a set of  unconstrained  problems. Moreover, if the steering control is scalar valued, then  the  optimal control can be expressed as a time-varying feedback law.\\
Let us consider the following system\\
\begin{equation}\label{P2-sys prin1}
    \left\{%
\begin{array}{ll}
\dot{y}(t)= Ay(t) +{\cal B} (u(t),y(t))   \\

    y(0)=y_0  \in X \\
\end{array}%
\right.
\end{equation}
		where
		\begin{itemize}
		\item{$A : D(A) \subset X  \mapsto X$ is the infinitesimal generator of a linear $C_0$- semi-group $S(t)$  on a real Hilbert space $X$ whose inner product and  corresponding norm are denoted respectively by $\langle.,.\rangle$ and $\Vert .\Vert$,}
\item   $u\in L^2(0,T;U)$, where $U$ is a real  Hilbert space equipped with inner product $\langle .,.\rangle_U$  and the corresponding norm $\Vert .\Vert_U$, and  $y$  is the corresponding mild solution to the control $u$,
		\item{$ {\cal B}: \ U\times X\rightarrow  X $  is a bounded  bilinear operator.}\end{itemize}
Let us now consider the  following  assumptions:\\
  $(a)$   For all $y\in X$ the mapping $u\mapsto {\cal B}(u,y)$ is compact,\\
 $(b)$ $A$ is the infinitesimal generator of a linear compact $C_0$- semigroup $S(t)$.\\
Note that assumption $(b)$ is systematically satisfied for $U=\mathbb{R}\cdot$\\
		The quadratic cost function $J$ to be minimized is defined by
		\begin{equation}\label{P2 cost obj}
			J(u)= \int_0^T \Vert y(t)\Vert^2 dt +\frac{r}{2} \int_0^T\Vert u(t)\Vert_U^2  dt\cdot
		\end{equation} Here, $r>0$ and $u$ belongs to the set of admissible control
$$U_{ad}=\{ u\in V\  \ /  \  \ y(T)=y_d\},$$
  where  $V$ is a closed convex subset of $L^2(0,T;U)$ and $y_d\in X$ is the desired state.\\
The optimal control problem may be stated as follows		
        \begin{equation*}\label{P2prob opt obj}
   (P)         \ \ \ \ \        \left\{%
\begin{array}{ll}
min J(u)  \\
    u\in U_{ad} \\
\end{array}%
\right.
\end{equation*}
     In order to solve the problem $(P)$, let us introduce the following auxiliary cost function\\
\begin{equation*}\label{P2 cost inter}
J_\epsilon(u)= \Vert y(T)-y_d\Vert^2 + \epsilon J(u),
 \end{equation*}
 where $\epsilon>0$,  and  let us consider the following optimal control problem
        \begin{equation*}\label{P2 prob opt inter}
   (P_\epsilon )\ \ \ \ \ \  \left\{%
\begin{array}{ll}
min J_\epsilon(u)  \\
    u\in V\\
\end{array}%
\right.
\end{equation*}
This paper is organized as follows: In Section 2, we will first provide a solution to the  auxiliary problem $(P_\epsilon)$.  This result is then applied to build a solution of the problem $(P)$. We will further  provide  sufficient  conditions  on  the  operators $A$ and $B$  under which the  solution of the problem $(P)$ can be expressed as a time-varying feedback law. Section 3 is devoted to examples and simulations.

\section{ Characterisation of the optimal control }
\subsection{Preliminary}
  Let us recall the notion of attainability.
  \begin{df}\
\begin{itemize}
\item A target state $y_d\in X$ is approximately attainable for the  system (\ref{P2-sys prin1}), if for all $\varepsilon> 0$ there exists $u_\varepsilon\in V$ such that $\Vert y_{u_\varepsilon}(T)-y_d\Vert \leq \varepsilon\cdot$
\item A target state $y_d\in X$ is exactly attainable for the  system (\ref{P2-sys prin1}), if there exists $u\in V$ such that $y_u(T)=y_d\cdot$
\end{itemize}

  \end{df}
  The following lemma provides a continuity property of the solution $y$ with respect to the control $u$.
  \begin{lem} \label{2}
  If one of the assumption $(a)$ or $(b)$ hold, then for any  sequence $(u_n)\subset L^2(0,T;U)$ such that $u_n\rightharpoonup u$ in $L^2(0,T;U)$,
  we have
  \begin{equation*}
  \lim_{n\rightarrow +\infty}\sup_{0\leq t \leq T}\Vert y_n(t)-y(t)\Vert=0,
  \end{equation*}
  where $y_n$ and $y$ are the mild solutions of the system (\ref{P2-sys prin1}) respectively corresponding to $u_n$ and $u\cdot$
  \end{lem}
\textbf{ Proof\\}
First, let us recall that  for all $u\in L^2(0,T;U)$, the system (\ref{P2-sys prin1}) has a unique mild solution corresponding to $u$, which is given by the following variation of constants formula (see e.g. \cite{Li and Young}, p. 66):\\
\begin{equation*}
y(t)=S(t)y_0+\int_0^tS(t-s){\cal B}(u(s),y(s))ds\cdot
\end{equation*}
Thus, the solutions  $y_n$ and $y$  of the system (\ref{P2-sys prin1}) respectively corresponding to $u_n$ and $ u$ satisfy the following formula for $t\in [0,T]$
\begin{equation*}
\begin{array}{r c l}
y_n(t)-y(t)&=&\int_0^tS(t-s)\bigg ({\cal B}(u_n(s),y_n(s))-{\cal B}(u(s),y(s))\bigg )ds\cdot
\end{array}
\end{equation*}
Then, for all $t\in [0,T]$ we have
\begin{equation*}
\Vert y_n(t)-y(t)\Vert \leq \Vert \int_0^t S(t-s){\cal B}(u_n(s)-u(s),y(s))ds\Vert + \Vert {\cal B}\Vert\int_0^t\Vert S(t-s)\Vert  \Vert u_n(s)\Vert_U  \Vert y_n(s)-y(s)\Vert ds
\end{equation*}
Applying  the Gronwall lemma (see Theorem 1 in \cite{gronwall}) yields
\begin{equation}\label{P2 y_n ver y^* 1}
\Vert y_n(t)-y(t)\Vert \leq \sup_{t\in [0,T]}\bigg (\Vert \int_0^t S(t-s){\cal B}(u_n(s)-u(s),y(s))ds\Vert \bigg)\exp\big (\Vert {\cal B}\Vert\int_0^t\Vert S(t-s)\Vert  \Vert u_n(s)\Vert_U  ds\big )
\end{equation}
Using the weak convergence  of  $u_n$ in $L^2(0,T;U)$ and the fact that the semi-group $S(t)$ is bounded on the entire finite interval $[0,T]$, we have for some $M>0$
\begin{equation}\label{P2 y_n ver y^* 2}
\exp\big (\Vert {\cal B}\Vert \int_0^t\Vert S(t-s)\Vert  \Vert u_n(s)\Vert_U ds\big )\leq M ,\ \  \ \forall t\in [0,T]\cdot
\end{equation}
\\
\textbf{$1^{st}$case : Assume that  $(a)$ holds.}\\
The weak convergence of $u_n\rightharpoonup u$ in $L^2(0,T;U)$ implies that ${\cal B}(u_n(.),y(.))$ strongly converge to $ {\cal B}(u(.),y(.))$ in $L^2(0,T;X)\cdot$\\
Then, we conclude that
\begin{equation}\label{P2 case1}
\lim_{n\rightarrow +\infty}\sup_{0\leq t\leq T}\Vert\int_0^tS(t-s){\cal B}(u_n(s)-u(s),y(s))ds\Vert=0\cdot
\end{equation}
It follows from (\ref{P2 y_n ver y^* 1}), (\ref{P2 y_n ver y^* 2}) and ( \ref{P2 case1})  that\\
\begin{equation*}
\lim_{n\rightarrow +\infty}\sup_{0\leq t \leq T}\Vert y_n(s)-y(s)\Vert=0\cdot
\end{equation*}
\textbf{$2^{nd}$ case : Assume that $(b)$ holds.}\\
According to Theorem 3.9 in \cite{Brezis}, the weak convergence : $u_n\rightharpoonup u$ in $L^2(0,T;U)$ implies the following weak convergence :  ${\cal B}(u_n(.),y(.))\rightharpoonup {\cal B}u(.),y(.))$ in $L^2(0,T;X)\cdot$\\
Moreover, the weak convergence of ${\cal B}(u_n(.),y(.))\rightharpoonup {\cal B}u(.),y(.))$ in $L^2(0,T;X)$ gives (see Corollary 3.3 of \cite{Li and Young}):
\begin{equation}\label{P2 y_n ver y^* 3}
\lim_{n\rightarrow +\infty}\sup_{0\leq t\leq T}\Vert\int_0^tS(t-s){\cal B}(u_n(s)-u(s),y(s))ds\Vert=0\cdot
\end{equation}
It follows from (\ref{P2 y_n ver y^* 1}), (\ref{P2 y_n ver y^* 2}) and ( \ref{P2 y_n ver y^* 3})  that\\
\begin{equation*}
\lim_{n\rightarrow +\infty}\sup_{0\leq t \leq T}\Vert y_n(s)-y(s)\Vert=0\cdot
\end{equation*}
\subsection{Optimal control for the problem $P_\epsilon$}

The following result discusses the existence of the optimal control related to  the auxiliary problem $(P_\epsilon)$.
\begin{thm} \label{thm3} \  \\
 Let one of the assumptions $(a)$ or $(b)$ hold.
\begin{itemize}
\item If  $V=\{u\in L^2(0,T;U)  /  \Vert u\Vert_U \leq M\} $ for some $M>0$, then there exists an optimal control for the problem $(P_\epsilon),$ which satisfies the following formula:
\end{itemize}

 \begin{equation*}\label{P2 optimal control t;x born}
u^*(t)=-\left( \frac{\Vert \epsilon r u^*(t)+({\cal B}(.,y^*(t))^*\phi(t)\Vert_U}{M}+\epsilon r\right)^{-1}{\cal B}(., y^*(t))^*\phi(t),
\end{equation*}
		 where  $\phi$ is the mild solution of the following adjoint system
        \begin{equation}\label{P2 adjoint pontr}
\begin{cases}
\dot{\phi}(t)=-A^*\phi(t)-{\cal B}^*(u^*(t),\phi(t))-2\epsilon y(t)\\
\phi(T)=2(y(T)-y_d)
\end{cases}
\end{equation}
 ${\cal B}^*(u^*(t),.)$ being the adjoint of the operator ${\cal B}(u^*(t),.)\cdot$
		
\begin{itemize}
\item If  $V=L^2(0,T;U)$, then the control defined by
\end{itemize}

		\begin{equation*}
		u^*(t)=-\frac{1}{\epsilon r} ({\cal B}(. , y^*(t))^*\phi(t)
		\end{equation*}
        is a solution of the problem $(P_\epsilon)$, where  $\phi$ is the mild solution of the adjoint system (\ref{P2 adjoint pontr}).\\
        \end{thm}
       \textbf{ Proof:}\\
       First let us show the existence of a solution of the problem $(P_\epsilon)$. \\
     Since the set  $\{J_\epsilon(u)/ u\in V\}\subset \mathbb{R}^+$ is not empty and bounded from below,  it admits a lower bound $J^*$.
		Let $ (u_n )_{n\in\mathbb{N}}$ be a minimizing sequence such that $J_\epsilon(u_n)\rightarrow J^*$. \\
		Then the sequence  $(u_n)$ is bounded, so it admits a sub-sequence still denoted by $(u_n)$, which weakly converges to   $u^*\in V$ .\\
		Let $y_n$ and $y^*$ be the solutions of (\ref{P2-sys prin1}) respectively corresponding to $u_n$ and $u^*$.\\
		From Lemma \ref{2} we have
		\begin{equation}\label{P2 limit ball}
		\lim_{n\rightarrow +\infty}\Vert y_n(t)-y^*(t)\Vert=0,  \   \forall t\in [0,T] \cdot
		\end{equation}
Since the norm  $\Vert . \Vert$ is lower semi-continuous, it follows from (\ref{P2 limit ball}) that for all $t\in [0,T]$\\
		\begin{equation*}
			\Vert y^*(t)\Vert^2= \lim_{n\rightarrow+\infty}inf\Vert y_n(t)\Vert^2\cdot
		\end{equation*}
		Applying Fatou's lemma we get \\
		\begin{equation}\label{P2 y inf Ja}
			\int_0^T\Vert y^*(t)\Vert^2dt= \lim_{n\mapsto +\infty}inf\int_0^T\Vert y_n(t)\Vert^2dt \cdot
		\end{equation}
		Since $R: u\mapsto \int_0^T\Vert u(t)\Vert_U^2dt $   is convex and lower semi-continuous with respect to weak topology, we have (see Corollary III.8 of \cite{Brezis})
		\begin{equation}\label{P2 u inf Ja}
			R(u^*)\leq \lim_{n\rightarrow+\infty}\inf R(u_n)\cdot
		\end{equation}Combining the  formulas (\ref{P2 limit ball}) , (\ref{P2 y inf Ja}) and (\ref{P2 u inf Ja}) we deduce that
		\begin{equation*}
      	\begin{aligned}
   J_\epsilon(u^*) & =  \Vert y^*(T)-y_d\Vert^2+\epsilon\int_0^T\Vert y(t)\Vert^2dt +\frac{\epsilon r}{2}\int_0^T\Vert u^*(t)\Vert_U^2dt\\
            &\leq  \lim_{n\rightarrow+\infty}\inf\Vert y_n(T)-y_d\Vert^2+\epsilon\lim_{n\rightarrow +\infty}inf \int_0^T\Vert y_n(t)\Vert^2dt+\frac{\epsilon r}{2}\lim_{n\rightarrow +\infty}inf \int_0^T\Vert u_n(t)\Vert_U^2dt\\\\
            &\leq  \lim_{n\rightarrow +\infty}infJ_\epsilon(u_n)\\\\
            &\leq J^*\cdot
            	\end{aligned}
		\end{equation*}
        We conclude that  $J_\epsilon(u^*)=J^*$ and so $u^*$ is a solution of the problem $(P_\epsilon)$.\\
    Let us proceed to the characterisation of the optimal control.
\begin{enumerate}
\item    \textbf{ The case  $V=\{ u\in L^2(0,T;U) \ \ /  \  \Vert u\Vert _{L^2(0,T,U)}\leq M\}\cdot$ }\
\end{enumerate}
Let $f_0: X\times U \mapsto \mathbb{R}$  be defined by
       \begin{equation*}
       f_0(y,u)=\epsilon \bigg (\Vert y\Vert^2+ \frac{r}{2}\Vert u\Vert_U^2\bigg ), \ \forall (y,u)\in X \times U\cdot
       \end{equation*}
       Then, the cost function $J_\epsilon $ takes the form
       \begin{equation*}
       J_\epsilon(u)=\Vert y(T)-y_d\Vert^2 +\int_0^T f_0(y(t),u(t))dt\cdot
       \end{equation*}
       Since $V$ is bounded, by application of Pontryagin's maximum principle (see Theorem 5.2 p. 258 in \cite{Li and Young} and Theorem 6.1 p. 162 in \cite{cannarsa 92} ), we find that for any  solution  $u^*$ of the problem $(P_\epsilon)$ there exists a  function $\phi$ solution of the following adjoint system
       \begin{equation*}
\begin{cases}
\dot{\phi}(t)=-A^*\phi(t)-{\cal B}^*(u^*(t),\phi(t))-2\epsilon y^*(t)\\
\phi(T)=2(y^*(T)-y_d)
\end{cases}
\end{equation*}
and satisfies the following condition
\begin{equation}\label{P2 principe max}
   H(t,u^*(t),y^*(t),\phi(t))=\min_{u\in V}H(t,u(t),y^*(t),\phi(t)),
   \end{equation}
where
\begin{equation*}
H(t,u(t),y^*(t),\phi(t))=f_0(u(t),y^*(t))+\langle \phi(t),{\cal B}(u(t),y^*(t))\rangle\cdot
\end{equation*}
By differentiating the function $u\mapsto H(u)=H(t,u(t),y^*(t),\phi(t))$, we have
\begin{equation*}
H'(u)(t) = \epsilon r u(t)+{\cal B}(. ,y^*(t))^*\phi(t),
\end{equation*}
where  $({\cal B}(. ,y^*(t))^*:X\mapsto U$ is the adjoint of the operator ${\cal B}(. ,y^*(t))$.\\
If   $\Vert u^*\Vert_{L^2(0,T;U)}<M$, then we conclude that
\begin{equation}\label{P2 H' 0}
u^*(t)=-\frac{1}{\epsilon r} {\cal B}(. , y^*(t))^*\phi(t)\cdot
\end{equation}
If   $\Vert u^*\Vert_{L^2(0,T;U)}=M$, we can distinguish two cases,  if $H'(u^*)=0$ then the control is given by (\ref{P2 H' 0})  and if  $H'(u^*)\neq 0$, then we proceed as follows:\\
Let
$v_1(t)=\frac{1}{M}u^*(t)$  and $v_2(t)=-\frac{1}{\Vert H'(u^*)\Vert_{L^2(0,T;U)}}H'(u^*)(t)$. We will show that $v_1=v_2\cdot$\\
For all $u\in V$ we have
\begin{equation*}
\langle v_1,u \rangle_{L^2(0,T;U)} \leq \Vert v_1\Vert_{L^2(0,T;U)} \Vert u\Vert_{L^2(0,T;U)} \leq M\ \ \ \ and \ \ \ \langle v_1,u^*\rangle_{L^2(0,T;U)}=M\cdot
\end{equation*}
So we conclude that
\begin{equation*}
\forall u\in V, \ \ \ \ \langle v_1,u\rangle_{L^2(0,T;U)}\leq \langle v_1,u^*\rangle_{L^2(0,T;U)}\cdot
\end{equation*}
Moreover, the fact that $V$ is convex, implies
\begin{equation*}
\forall u\in V, \ \  \forall \lambda \in [0,1],  \ \  u^* + \lambda(u-u^*) \in V\cdot
\end{equation*}
Then since $u^*$ is a solution of the problem $(P_\epsilon)$, we  derive from
%find  by the  maximum principle given by
(\ref{P2 principe max})
% (as $\lambda \rightarrow 0^+$)\\
\begin{equation}\label{P2 H' VN}
\begin{array}{ccc}
 H(u^*)&\leq& H(u^*+\lambda (u-u^*)) \\
  &\leq& H(u^*)+\langle H'(u^*),\lambda (u-u^*)\rangle_{L^2(0,T;U)} \\
   &+ &\lambda \Vert u^*-u\Vert_{L^2(0,T;U)}\theta(\lambda \Vert u^*-u\Vert_{L^2(0,T;U)}),  \ \forall \lambda\in [0,1], \ \forall u\in V
\end{array}
\end{equation}
where the function $\theta$ is such that
\begin{equation}\label{P2 H' teta}
\lim_{\lambda \rightarrow 0^+}\theta(\lambda \Vert u^*-u\Vert_{L^2(0,T;U)} )=0\cdot
\end{equation}
From (\ref{P2 H' VN}) and (\ref{P2 H' teta}) it comes
\begin{equation*}
\langle H'(u^*),u\rangle_{L^2(0,T;U)} \ge \langle H'(u^*), u^*\rangle_{L^2(0,T;U)}\cdot
\end{equation*}
So, we conclude that
\begin{equation*}
\forall u\in V_{ad},  \ \ \ \ \langle v_2, u\rangle_{L^2(0,T;U)}\leq \langle v_2,u^*\rangle_{L^2(0,T;U)}\cdot
\end{equation*}
Taking into account that $\sup_{u\in V}\langle v_2,u\rangle_{L^2(0,T;U)}=M$, we deduce that  $\langle v_2,u^*\rangle_{L^2(0,T;U)}=M$ and that
\begin{equation*}
\langle \frac{1}{2}(v_1+v_2),u^*\rangle_{L^2(0,T;U)}=\frac{1}{2}\langle v_1,u^*\rangle_{L^2(0,T;U)} + \frac{1}{2}\langle v_2,u^*\rangle_{L^2(0,T;U)}=M,
\end{equation*}
then
\begin{equation*}
\Vert \frac{1}{2}(v_1+v_2)\Vert_{L^2(0,T;U)}\ge 1\cdot
\end{equation*}
It follows that
$$\Vert (v_1+v_2)\Vert_{L^2(0,T;U)}=\Vert v_1\Vert_{L^2(0,T;U)} + \Vert v_2\Vert_{L^2(0,T;U)} $$
and that   $v_1=v_2$.\\
Furthermore, we have
\begin{equation}\label{P2 H' no}
\frac{1}{M}u^*(t)=-\frac{1}{\Vert H'(u^*)\Vert_{L^2(0,T;U)}}H'(u^*)(t)\cdot
\end{equation}
According to (\ref{P2 H' 0}) and (\ref{P2 H' no}) we have
\begin{equation*}
u^*(t)=\frac{-1}{\frac{\Vert H'(u^*)\Vert_{L^2(0,T;U)}}{M}+\epsilon r} {\cal B}(. , y^*(t))^*\phi(t),
\end{equation*}
where
\begin{equation*}
H'(u)(t) = \epsilon r u(t)+{\cal B}(. , y^*(t))^*\phi(t)\cdot
\end{equation*}

\textbf{ 2.    The case $V=L^2(0,T;U)$}.\\\\
From the first part of the proof, there exists a  solution $u^*$ of the problem $(P_\epsilon)$.\\
Let us consider the closed convex space $$V^*=\{ u\in L^2(0,T;U) \   / \ \Vert u\Vert_{L^2(0,T;U)} \leq \Vert u^*\Vert_{L^2(0,T;U)} +1\}\cdot$$
It is clear that $u^*(t)\in \mathring{V}^*$, then from the first case, we have $H'(u^*)=0$, which leads to
\begin{equation*}
		u^*(t)=-\frac{1}{\epsilon r} {\cal B}(.,y^*(t))^*\phi(t),
		\end{equation*}
        where  $\phi$ is the mild solution of the adjoint system (\ref{P2 adjoint pontr}).\\
This achieves the proof of Theorem \ref{thm3}.
\subsection{Sequential characterization of the solution of the problem $(P)$}
        In the sequel, we take a decreasing  sequence $(\epsilon_n)$ such that $\epsilon_n\rightarrow 0$ with corresponding sequence of controls  $(u_n^*)$  solutions of problems $(P_{\epsilon_n})$.
\begin{thm}\label{thm4}
Assume  that $V$ is bounded and let $y_d$  be an approximately attainable state by a control from $V$. Then the problem $(P)$ posses a solution. Moreover any weak limit value of $(u_n^*)$ in $L^2(0,T,U)$ is a solution of   $(P)$.

\end{thm}
\textbf{Proof:}\\
Since  $V$  is bounded, we deduce that the sequence  $(u_n^*)$ is bounded, so it admits a  weakly converging subsequence, denoted by $(u_n^*)$ as well. Let $u^*$ be a weak limit value of $(u_n^*)$ in $V$.\\
The remainder of the proof is divided into three steps\\
\textbf{Step 1: $y_d$ is exactly attainable $U_{ad}\ne \emptyset \cdot$}\\
Let us consider the following problem \\
\begin{equation}\label{P2 prob non vide}
    \left\{%
\begin{array}{ll}
min \ \  \Vert y_u(T)-y_d\Vert^2 \\
    u\in V\\
\end{array}%
\right.
\end{equation}
The set  $\{\Vert y_u(T)-y_d\Vert^2 / u\in V\}\subset \mathbb{R}^+$ is not empty and bounded from below, so it admits a lower bound $J_d$.\\
		Let $ (v_n )_{n\in\mathbb{N}}$ be a minimizing sequence such that $\Vert y_{v_n}(T)-y_d\Vert^2 \underset{n\rightarrow+\infty}{\longrightarrow}  J_d$. \\
		Since  $V$ is bounded, we deduce that the sequence  $(v_n)$ is bounded, so it admits a weakly converging subsequence to  $v\in V$ still denoted by $(v_n)$.\\
By Lemma \ref{2}, we have for all $t\in [0,T]$
\begin{equation*}
\lim_{n\rightarrow +\infty}\Vert y_{v_n}(t)-y_v(t)\Vert = 0
\end{equation*}
then, we conclude that
\begin{equation}\label{P2 min bor}
			\Vert y_v(T)-y_d\Vert^2=\lim_{n\rightarrow + \infty} \Vert y_{v_n}(T)-y_d\Vert^2 =J_d= \min_{u\in V}\Vert y_u(T)-y_d\Vert^2
		\end{equation}
So the control  $v$ is a solution of the problem (\ref{P2 prob non vide}).\\
Since the system (\ref{P2-sys prin1}) is approximately attainable, we have
\begin{equation}\label{P2 approxi def}
\forall \varepsilon>0 ,\ \ \exists v_{\varepsilon}\in V \ \  / \  \  \   \Vert y_{v_\epsilon}(T)-y_d\Vert\leq \varepsilon
\end{equation}
According to (\ref{P2 min bor})  and (\ref{P2 approxi def}), we get
\begin{equation*}
\forall \epsilon > 0, \ \ \exists v_{\varepsilon}\in  V,   \ \ \  \Vert y_v(T)-y_d\Vert\leq \Vert y_{v_\epsilon}(T)-y_d\Vert \leq \epsilon
\end{equation*} So we conclude that $\Vert y_v(T)-y_d\Vert =0$ and hence  $v\in U_{ad}$. \\

\textbf{Step 2:  $  \forall v\in U_{ad},\ \ J(u^*)\leq J(v)\cdot$}\\
Taking into account that  $u_n^*$ is a solution of the problem $(P_{\epsilon_n})$ and $y_n^*$ is the corresponding solution of the system (\ref{P2-sys prin1}), we get for all $v\in U_{ad}$
\begin{equation*}
J_{\epsilon_n}(u_n^*)=\Vert y_n^*(T)-y_d\Vert^2 + \epsilon_nJ(u_n^*)\leq J_{\epsilon_n}(v)
\end{equation*}from which, it comes
\begin{equation*}
\begin{array}{r c l}
				\epsilon_nJ(u_n^*)  &\leq & J_{\epsilon_n}(v)-\Vert y_n^*(T)-y_d\Vert^2 \\\\
				&\leq & \epsilon_nJ(v)
			\end{array}
\end{equation*}So we find
\begin{equation}\label{P2 contro inf 2}
J(u_n^*)\leq J(v)   \ for \ all \ \ v\in U_{ad}\cdot
\end{equation}
Let  $y^*$ be the solution of system  (\ref{P2-sys prin1})  corresponding to $u^*$.\\
		Since  $u_n\rightharpoonup u^*$ in $L^2(0,T;U)$, we have by Lemma \ref{2}
		\begin{equation}\label{P2 yn tend y*}
		\lim_{n\rightarrow +\infty}\Vert y_n^*(t)-y^*(t)\Vert=0,  \   \forall t\in [0,T] \cdot
		\end{equation}The norm  $\Vert . \Vert$ is lower semi-continuous, it follows that for all $t\ge 0$ we have\\
		\begin{equation*}
			\Vert y^*(t)\Vert^2=\lim_{n\rightarrow+\infty}inf\Vert y^*_n(t)\Vert^2\cdot
		\end{equation*}
		Applying Fatou's lemma we get\\
		\begin{equation}\label{P2 y inf Jco}
			\int_0^T\Vert y^*(t)\Vert^2dt=\lim_{n\rightarrow +\infty}inf\int_0^T\Vert y^*_n(t)\Vert^2dt \cdot
		\end{equation}
		The function  $R$  is lower semi-continuous and convex, it follows from \cite{Brezis} that
		\begin{equation}\label{P2 u inf Jco}
			R(u^*)\leq \lim_{n\rightarrow+\infty}\inf R(u_n^*)\cdot
		\end{equation}By  the  inequalities (\ref{P2 y inf Jco} ) and (\ref{P2 u inf Jco}) we deduce that
		\begin{equation}\label{P2 limi inf born 2}
			J(u^*)\leq \lim_{n\rightarrow +\infty}inf(J(u_n))\cdot
		\end{equation}Combining  ( \ref{P2 contro inf 2}) and (\ref{P2 limi inf born 2}) we deduce that
        \begin{equation*}
        J(u^*)\leq J(v)\cdot
        \end{equation*}\textbf{Step 3 :  $u^*\in U_{ad}\cdot$}\\
  According to the inequality (\ref{P2 contro inf 2}), we deduce that $J(u_n^*)$ is bounded and
  \begin{equation*}
  \lim_{n\rightarrow +\infty} \Vert y^*_n(T)-y_d\Vert^2= \lim_{n\rightarrow +\infty} J_{\epsilon_n}(u_n^*)\leq \lim_{n\rightarrow +\infty}J_{\epsilon_n}(v)=\Vert y_v(T)-y_d\Vert^2=0\cdot
  \end{equation*}
Then, taking into account the formula (\ref{P2 yn tend y*}), we derive via the continuity of the norm that
\begin{equation*}
\lim_{n\rightarrow +\infty} \Vert y_n^*(T)-y_d\Vert = \Vert y^*(T)-y_d\Vert \leq \Vert y_v(T)-y_d\Vert = 0 \cdot
\end{equation*}
Consequently, $y^*(T)=y_d$ and the control $u^*$ is a solution of problem $(P)$.\\
\begin{thm} \label{thm5}
If $U_{ad}\ne \emptyset$ , then there exists a  solution  $u^*$ of the problem $(P)$. Furthermore, any weak limit value of the solution $(u_n^*)$ of $(P_{\epsilon _n})$ in $L^2(0,T;U)$ is a solution of $(P)$.
\end{thm}
\textbf{Proof:}\\
Let $v\in U_{ad}$ . Then keeping in mind that $u_n^*$ is the solution of the problem $(P_{\epsilon_n})$ corresponding to $\epsilon_n$, we can see that
\begin{equation*}
J_{\epsilon_n}(u_n^*)\leq J_{\epsilon_n}(v)= \epsilon_nJ(v)
\end{equation*}
It follows that
\begin{equation*}
\epsilon_nJ(u_n^*)=J_{\epsilon_n}(u_n^*)-\Vert y_n^*(T)-y_d\Vert^2\leq J_{\epsilon_n}(u_n^*)\leq \epsilon_nJ(v)
\end{equation*}
Using the definition of the cost $J$ given by (\ref{P2 cost obj}), the last equality gives
\begin{equation}\label{contro inf 2 }
r\int_0^T\Vert u_n^*(t)\Vert_U ^2dt\leq J(u_n^*)\leq J(v)\cdot
\end{equation}
We deduce that the sequence  $(u_n^*)$ is bounded, so it admits a  weakly converging subsequence in $V$, also denoted by $(u_n^*)$.
Let $u^*$ be a weak limit value of $(u_n^*)$ in $V$ and let $y^*$ be the solution of system (\ref{P2-sys prin1})  corresponding to  $u^*$.\\
		Since  $u_n\rightharpoonup u^*$ in $L^2(0,T;U)$, we have by Lemma \ref{2}
		\begin{equation*}\label{yn tend y*}
		\lim_{n\rightarrow +\infty}\Vert y_n^*(t)-y^*(t)\Vert=0,  \   \forall t\in [0,T] \cdot
		\end{equation*}
		Similarly to the  proof of Theorem \ref{thm4}  we can show that
        \begin{equation*}
        J(u^*)\leq J(v)\cdot
        \end{equation*}
  According to the inequality (\ref{contro inf 2 }), we deduce that $J(u_n^*)$ is bounded and
  \begin{equation*}
  \lim_{n\rightarrow +\infty}J_{\epsilon_n}(u_n^*)=\lim_{n\rightarrow +\infty} \Vert y_n^*(T)-y_d\Vert^2 \leq \Vert y_v(T)-y_d\Vert^2\cdot
  \end{equation*}
Hence
\begin{equation*}
\lim_{n\rightarrow +\infty} \Vert y_n^*(T)-y_d\Vert = \Vert y^*(T)-y_d\Vert \leq \Vert y_v(T)-y_d\Vert = 0
\end{equation*}
We conclude that $u^*\in U_{ad}\cdot$\\

\subsection{Optimal feedback control }
In this part we will try to express the optimal control $u^*$ of the problem $(P)$ as a time-varying feedback law for the class of commutative bilinear systems with scalar control \cite{Alami,Wei}.\\
Assume that $U=\mathbb{R}$, then we can write the system (\ref{P2-sys prin1}) as follows
\begin{equation*}
\begin{cases}
\dot{y}(t)= Ay(t) + u(t)By(t)   \\

    y(0)=y_0  \in X \\
\end{cases}
\end{equation*}
where {$A : D(A) \subset X  \mapsto X$ is the infinitesimal generator of a linear  $C_0$- semi-group $S(t)$,}  $B$ is a bounded linear operator and $u\in V:=L^2(0,T)$ .\\

\begin{thm} \label{thm6}
Assume that $A$ and $B$ commute with each other and that $ U_{ad}\ne \emptyset$. Let $v\in U_{ad} $ and let $y_0\in X$ be such that $S(T)y_0\not\in Ker(B)$. Then for any solution  $u^*$ of the problem $(P)$, we have  the following formula
\begin{equation*}
u^*(t)= \frac{1}{T}\int_0^Tv(s)ds+\frac{2}{Tr}\int_0^T\int_\alpha^T\langle y^*(s),By^*(s)\rangle ds\  d\alpha-\frac{2}{r}\int_t^T\langle y^*(s),By^*(s)\rangle ds
\end{equation*}
\end{thm}
\textbf{Proof:}\\
Let us consider the system (\ref{P2-sys prin1}) in the  time horizon $[0,T]$, and  let  $A_k=kA(kI-A)^{-1}$ be the Yosida approximation of the operator $A$.
		Let $y_k$ and $\phi_k$  be the respective solutions to (\ref{P2-sys prin1}) and (\ref{P2 adjoint pontr}) with $A_k$ instead of $A$. For $u\in L^2(0,T)$,
		since $A_k$ is bounded, we have $ y_k, \phi_k \in
		H^1(0, T)$ and
		\begin{equation*}
			\begin{array}{r c l}
				\langle\dot{\phi}_k(t),By_k(t)\rangle+\langle \phi_k(t),B\dot{y}_k(t)\rangle &=&\langle -A_k^*\phi_k(t)-u(t)B^*\phi_k(t)-2\epsilon
                y_k(t),By_k(t)\rangle \\\\ &+&\langle B^*\phi_k(t), A_ky_k(t)+u(t)By_k(t)\rangle\\\\
				&=&\langle\phi_k(t),BA_ky_k(t)-A_kBy_k(t)\rangle-2\epsilon\langle y_k(t),By_k(t)\rangle\cdot\\
			\end{array}
		\end{equation*}
		Thus
		\begin{equation}\label{P2 differ phi commut}
			\langle\dot{\phi}_k(t),By_k(t)\rangle+\langle \phi_k(t),B\dot{y}_k(t)\rangle=\langle\phi_k(t),[B,A_k]y_k(t)\rangle-2\epsilon\langle y_k(t),By_k(t)\rangle
		\end{equation}
where $[B,A_k]:=B A_k-  A_k B$.\\
		Integrating (\ref{P2 differ phi commut})  over $[t, T],$ we get
		\begin{equation*}\label{P1 phi By}
			\langle\phi_k(t),By_k(t)\rangle=2\langle y_k(T)-y_d,By_k(T)\rangle -\int_t^T\bigg (\langle\phi_k(s),[B,A_k]y_k(s)\rangle-2\epsilon\langle y_k(s),By_k(s)\rangle\bigg ) ds
		\end{equation*}
		Since $\phi_k\longrightarrow \phi$ and   $y_k\longrightarrow y$ strongly, we obtain  by letting $k\rightarrow +\infty$
    \begin{equation*}\label{P1 phi By lim}
			\langle\phi(t),By(t)\rangle=2\langle y(T)-y_d,By(T)\rangle+2\epsilon\int_t^T\langle y(s),By(s)\rangle ds\cdot
		\end{equation*}So, by Theorem \ref{thm3}, we conclude that the solution of the problem $(P_{\epsilon_n})$ corresponding to $\epsilon_n$, is given by
        \begin{equation}\label{P2 u_n feedback}
        u_n^*(t)=-\frac{1}{\epsilon_nr}\langle \phi_n(t),By_n^*(t)\rangle= -\frac{2}{\epsilon_nr}\langle y_n^*(T)-y_d,By_n^*(T)\rangle-\frac{2}{r}\int_t^T\langle y_n^*(s),By_n^*(s)\rangle ds\cdot
        \end{equation}
Let $v\in U_{ad}$. By Theorem \ref{thm5}, any  limit value $u^*$ of $u_n^*$  in $L^2(0,T)$ is a solution of the problem $(P)$.\\
Since $A$ and $B$ commute, we have the following formulas
\begin{equation*}
y_v(t)=S(t)\exp(B\int_0^tv(s)ds)y_0
\end{equation*}
and
\begin{equation*}
y^*(t)=S(t)\exp(B\int_0^tu^*(s)ds)y_0\cdot
\end{equation*}
Using the fact that $v,u^*\in U_{ad}$  and $\lim_{n\rightarrow +\infty } y_n^*(T)=y_d$ , we obtain
\begin{equation*}
\lim_{n\rightarrow + \infty} y_n^*(T) = y_u^*(T)=y_v(T)=y_d\cdot
\end{equation*}
Hence
\begin{equation*}
\lim_{n\rightarrow +\infty}S(T)\exp(B\int_0^Tu_n^*(t)dt)y_0=S(T)\exp(B\int_0^Tv(t)dt)y_0=S(T)\exp(B\int_0^Tu^*(t)dt)y_0\cdot
\end{equation*}
From  the assumption $S(T)y_0\not\in Ker (B)$, we deduce from the last inequalities that
\begin{equation*}
\lim_{n\rightarrow +\infty}\int_0^Tu_n^*(t)dt=\int_0^Tv(t)dt=\int_0^Tu^*(t)dt\cdot
\end{equation*}
Moreover, we deduce from the formula  (\ref{P2 u_n feedback}), that
\begin{equation*}
\begin{aligned}
\lim_{n\rightarrow +\infty}\ \ \int_0^Tu_n^*(t)dt &= \lim_{n\rightarrow +\infty} \int_0^T\bigg ( -\frac{2}{\epsilon_nr}\langle y_n^*(T)-y_d,By_n^*(T)\rangle-\frac{2}{r}\int_t^T\langle y_n^*(s),By_n^*(s)\rangle ds\bigg ) dt \\\\
&= \lim_{n\rightarrow +\infty} -\frac{2T}{\epsilon_nr}\langle y_n^*(T)-y_d,By_n^*(T)\rangle-\frac{2}{r}\int_0^T\int_t^T\langle y^*(s),By^*(s)\rangle ds dt
\end{aligned}
\end{equation*}
from which, we derive
\begin{equation}\label{P2 terme feedb}
\lim_{n\rightarrow +\infty} -\frac{2T}{\epsilon_nr}\langle y_n^*(T)-y_d,By_n^*(T)\rangle=\int_0^Tv(t)dt+ \frac{2}{r}\int_0^T\int_t^T\langle y^*(s),By^*(s)\rangle ds dt\cdot
\end{equation}
By (\ref{P2 u_n feedback}) and (\ref{P2 terme feedb}) we deduce that  $u_n^*(t)\rightarrow u^*(t)$ for all $t\in [0,T]$ and
\begin{equation*}
\begin{aligned}
\lim_{n\rightarrow +\infty}u_n^*(t)&=\lim_{n\rightarrow +\infty} -\frac{2}{\epsilon_nr}\langle y_n^*(T)-y_d,By_n^*(T)\rangle-\frac{2}{r}\int_t^T\langle y_n^*(s),By_n^*(s)\rangle ds\\\\
&= \frac{1}{T}\int_0^Tv(s)ds+\frac{2}{Tr}\int_0^T\int_\alpha^T\langle y^*(s),By^*(s)\rangle ds d\alpha-\frac{2}{r}\int_t^T\langle y^*(s),By^*(s)\rangle ds\\\\ &= u^*(t)\cdot
\end{aligned}
\end{equation*}
We conclude that
\begin{equation*}
u^*(t)= \frac{1}{T}\int_0^Tv(s)ds+\frac{2}{Tr}\int_0^T\int_\alpha^T\langle y^*(s),By^*(s)\rangle ds d\alpha-\frac{2}{r}\int_t^T\langle y^*(s),By^*(s)\rangle ds\cdot
\end{equation*}
\begin{rem} In the case where $S(t_1)$ is one to one for some $t_1>0$ and $y_0\not\in Ker(B), $ the assumption $S(T)y_0\not\in Ker(B)$ in Theorem 6 is satisfied.
\end{rem}

\section{Examples}
\subsection{Wave equation}
Let us consider the following wave equation
\begin{equation*}
\left\{
\begin{array}{lllll}
\frac{\partial^2}{\partial t^2}z(t,x) &=& \Delta z(t,x) +u (t,x)z(t,x), & \ \ t\in[0,T] \ and \ \ & x\in \Omega=(0,1)\\
z(t,0)&=&z(t,1)=0, & t\in [0,T] &\\
z(0,x)&=&z_0(x),\ &    &   x\in \Omega
\end{array}
\right.
\end{equation*}
where
\begin{itemize}
\item $u\in L^2(0,T,L^2(\Omega))$,
\item $T>4\max_{x\in \Omega} \vert x-x_0\vert$ for some $ x_0\in \mathbb{R}\setminus [0,1]$,
\item the desired state $z_d \in H_0^1(\Omega)\cap H^2(\Omega)$ is such that  $\frac{\Delta z_d}{z_d}\mathbb{1}_{(z_d\ne0)}\in L^\infty(\Omega)$, where $\mathbb{1}_{(z_d\ne0)} $ indicates the characteristic function of the set  $(z_d\ne0):=\{ x\in \Omega \ /\ \ z_d(x)\ne 0\}.$

\end{itemize}
    This system has the form of the system (\ref{P2-sys prin1}) if we take $y(t)=(z(t), \dot{z}(t))$, $X=H_0^1(\Omega)\times L^2(\Omega)$ with
$\langle (y_1,z_1),(y_2,z_2)\rangle_X=\langle y_1,y_2\rangle_{H_0^1(\Omega)} + \langle z_1,z_2\rangle_{L^2(\Omega)}$
and
\begin{equation*}
A=\begin{pmatrix}
0 & I\\
\Delta & 0
\end{pmatrix} \ \mbox{with} \ D(A)= H_0^1(\Omega)\cap H^2(\Omega)\times H_0^1(\Omega) \ \mbox{and}
\ \ B=\begin{pmatrix}
0 & 0\\
I & 0
\end{pmatrix}\cdot
\end{equation*}
Here $B$ is a compact linear bounded operator on $X$  and $A$ is the infinitesimal generator of a linear $C_0$- semi-group $S(t)$ of isometries (see \cite{ball 79},  p.176).\\
The quadratic cost function is given by
\begin{equation*}\label{100}
J(u)= \int_0^T (\Vert z(t)\Vert_{H_0^1(\Omega)}^2+\Vert \dot{z}(t)\Vert_{L^2(\Omega)}^2) dt + \frac{r}{2}\int_0^T\Vert u(t)\Vert_{L^2(\Omega)}^2 dt,
\end{equation*}
where $u(t):=u(t,\cdot)$ and $z(t):=z(t,\cdot).$ \\
According to  \cite{ouzahra2019wave}, there exists a control $v\in L^2(0,T;L^2(\Omega)$ such that the corresponding solution $z_v$ of the system (\ref{P2-sys prin1}) verifies $z_v(T)=z_d$.
Then,  according to Theorem \ref{thm5} there exists a control $u^*\in L^2(0,T,\mathbb{R})$,  which guarantees  the exact attainability of  $z_d$  at time $T$, and  is a solution of the problem $(P)$ with $U_{ad}=\{u\in L^2(0,T,L^2(\Omega))\  /   z(T)=z_d\}$.\\

\begin{rem} The optimal control of the bilinear wave equation   has been considered  in \cite{linag,Boukhari 2017} in the context of unconstrained endpoint.
\end{rem}

\subsection{Heat equation}
In this part  we study the optimal exact attainability  for the reaction-diffusion equation.\\
Let us consider  the following system

\begin{equation}\label{P2 chal}
            \left\{%
\begin{array}{ll}
					\frac{\partial}{\partial t}y(t,x)=\Delta y(t,x) +u(t,x)y(t,x), \ \ \ \ & in \ Q=\Omega \times(0,T),   \  T>0\\
                    y(t,0)=y(t,1)=0, \ \ \ \ \ \ \  \  \  \ \ \ \ \ \  \ \ \ \ \ \ \ \ \  \ & on \  (0,T)\\
                    y(0)=y_0 \ \ \  \ \ \ \ \ \ \ \ \ \ \ \ \ \  \ \ &  in \ \ \Omega
				\end{array}%
\right.
			\end{equation}
where $\Omega=(0,1) $ and $u\in L^2(0,T,U)$ is a control function.\\

\textbf{Case 1: Distributed control $(U=L^2(\Omega))$}\\

 Assume that $y_0, y_d \in L^2(\Omega)$ are such that
\begin{itemize}
\item for a.e. $x\in \Omega$ , $y_dy_0 \ge 0$, \\
\item for a.e. $x\in \Omega , \  y_0(x) = 0 \iff  y_d(x) = 0,$\\
\item  $a :=ln(\frac{y_d}{y_0}) {\bf 1}_{(y_0\ne0)} \in L^{\infty }(\Omega),$ where $\mathbb{1}_{(y_0\ne0)} $ indicates the characteristic function of the set  $(y_0\ne0):=\{ x\in \Omega \ /\ \ y_0(x)\ne 0\}.$

\item $ \frac{\Delta y_d}{y_d} 1_{{(y_d\ne0)}} \in L^\infty(\Omega), $\\
\item $\vert y_d\vert >0$  a.e. on some nonempty open subset $O$ of $\Omega$.
\end{itemize}
According to Theorem 2 in \cite{contr 2016}, there is a time $T$  for which $y_d$  is exactly attainable for the system (\ref{P2 chal}) using a control $v\in L^2(0,T, L^2(\Omega))$, so $U_{ad}\ne \emptyset $.
Then, according to Theorem \ref{thm5}, there exists a control $u^*$ which guarantees  the exact attainability of  $y_d$  at time $T$, and is  solution of the following problem
\begin{equation}\label{P2 Exemple}
    \left\{%
\begin{array}{ll}
min J(u)  \\
    u\in U_{ad}=\{u\in L^2(0,T,L^2(\Omega))\ \ \ / \ \ y_u(T)=y_d\} \\
\end{array}%
\right.
\end{equation}
More precisely  any weak limit of $u_n^*$ given by Theorem \ref{thm3} corresponding to sequence $(\epsilon_n)$ gives  a optimal control $u^*$ for (\ref{P2 Exemple}).\\

\textbf{Case 2: Scalar control $(U=\mathbb{R})$}\\
Here, we have $u(t,x)=u(t)\in \mathbb{R}$. \\
Assume that $y_0, y_d \in L^2(\Omega)$ are such that $y_d=\lambda y_0$ with $\lambda > 1$ and $y_0>0$, a.e in $\Omega$.
According to Theorem II 4 and Remark 4 in \cite{contr 2020}, there is a time $T$  for which $y_d$  is exactly attainable for the system (\ref{P2 chal}) using the control $v(t)=\frac{\lambda-1}{T+(\lambda-1)t}\in L^2(0,T, \mathbb{R})$, so $U_{ad}\ne \emptyset$. \\
By Theorem \ref{thm6}, there exists a feedback control $u^*\in L^2(0,T,\mathbb{R})$ which guarantees  the exact attainability of  $y_d$  at time $T$, and is  solution of the problem $(P)$ with $U_{ad}=\{u\in L^2(0,T,\mathbb{R})\ \ \ / \ \ y^*(T)=y_d\}$, and satisfies the following formula
\begin{equation*}
u^*(t)=\frac{1}{T}\ln(\lambda)+\frac{2}{Tr}\int_0^T\int_\alpha^T\Vert y^*(s)\Vert^2 ds\  d\alpha-\frac{2}{r}\int_t^T\Vert y^*(s)\Vert^2 ds\cdot
\end{equation*}
\subsection{Transport equation}
Let us consider the following transport problem
        \begin{equation}\label{P1 trans}
            \left\{%
\begin{array}{ll}
					\frac{\partial}{\partial t}y(t,x)= -\frac{\partial }{\partial x}y(t,x) +u (t)y(t,x), & t\in (0,T),\;  x\in \Omega=(0,+\infty)\\
					y(t,0)=0, & t\in (0,T) \\
                    y(0,x)=y_0(x), & x\in \Omega
				\end{array}%
\right.
			\end{equation}
		where $u\in L^2(0,T)$. Here the operator $A=-\frac{\partial}{\partial x}$ with the domain $D(A)=H_0^1(\Omega)$ generates a $C_0-$semi-group  of isometries $S(t)$ in $X=L^2(\Omega)$.
        Below, we will develop numerical  simulation for the example (\ref{P1 trans}). For this end, we take   $r=2$ , $T=9$, $y_0=x\exp(-x)$ and
        $$
        y_d(x)=
        \begin{cases}
        0, \ if \  x\leq 9\\
        (x-9)\exp(9-x),\ if \ x\ge 9
        \end{cases}
        $$
        then the control $v=0\in U_{ad}=\{u\in L^2(0,T)\ \ / \ \ y^*(T)=y_d\}\cdot$
      By Theorem \ref{thm6}, there exists a feedback control $u^*\in L^2(0,T)$ which guarantees  the exact attainability of  $y_d$  at time $T$. Moreover $u^*$ is  the solution of the problem $(P)$   and satisfies the following formula
\begin{equation}\label{P2 transport 1}
u^*(t)=\frac{1}{T}\int_0^T\int_\alpha^T\Vert y^*(s)\Vert^2 ds\  d\alpha-\int_t^T\Vert y^*(s)\Vert^2 ds\cdot
\end{equation}
In the Figure 1, we compare numerically the two controls $u^*$ and $v=0$ in term of the state at the finite time $T=9$. Moreover, we find $J(u^*)=1.2442$ and $J(v)=2.25\approx 2J(u^*)$.\\
\begin{figure}[h]
\centering
\includegraphics[width=15cm, height=6cm]{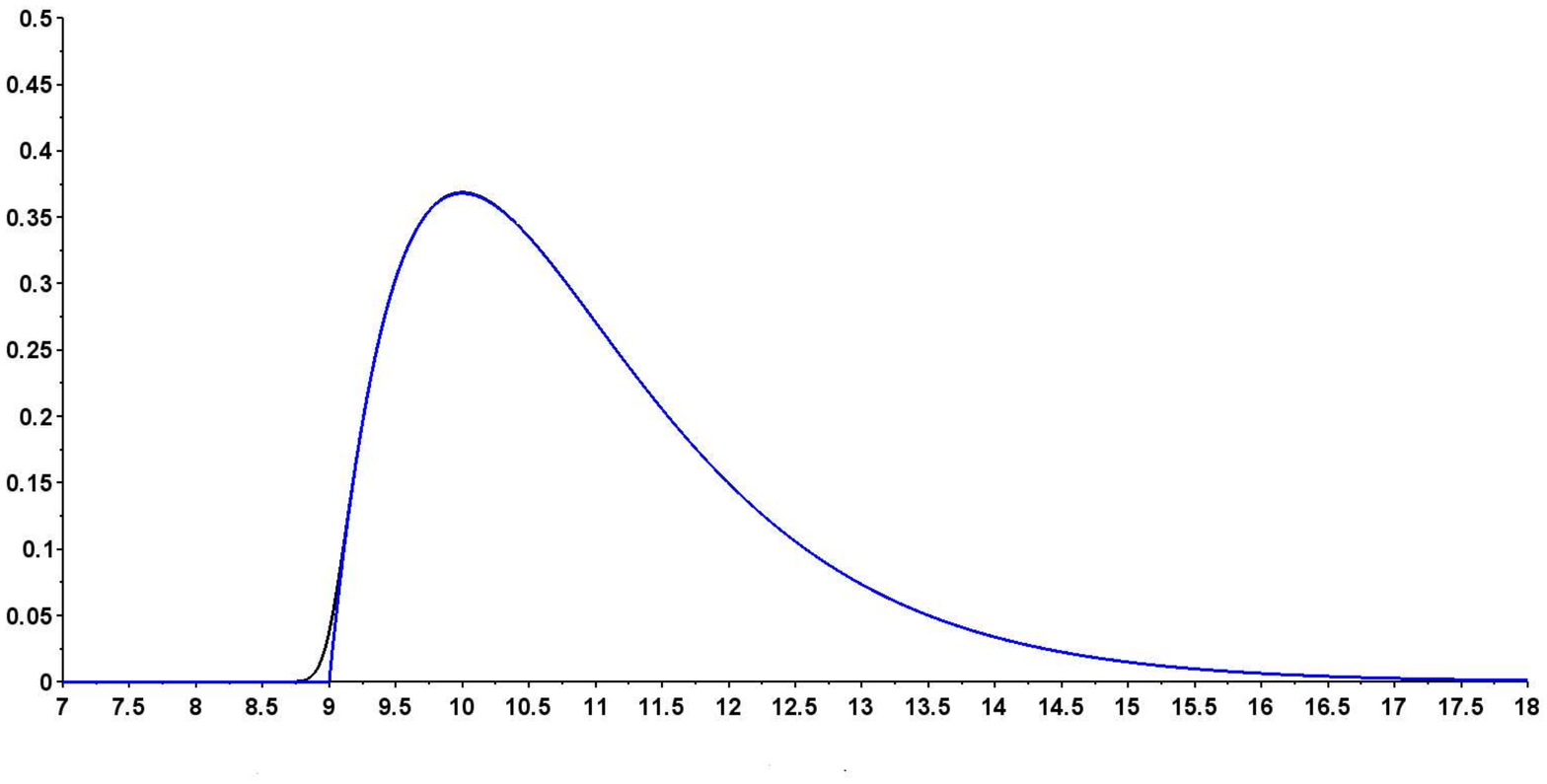}
\caption{The state $y_{u^*}(T)$  (black line), and the desired state $y_d$ (blue line) }
\end{figure}
We observe that the desired state is exactly attainable either  by using  the optimal control $u^*$ or the  control $v=0$.
However, the  control $ u^* $ leads to a lower cost than the zero control.\\

\begin{rem}Unlike the case of linear systems, the uniqueness of the   optimal control of the quadratic cost (\ref{P2 cost obj}) is not guaranteed in general when dealing with bilinear systems, which is due to the lack of convexity of the state w.r.t control. For instance, if we assume that  $r=0$ and   that $A=B $ is  a skew-adjoint matrix, we can see that the cost function is constant so we have an infinity of optimal controls. However,
in the case of the quadratic cost function  $J(u)=\int_0^T u^2(t)dt$, the uniqueness of the  optimal control is assured by the strict convexity of the cost $J$ (see \cite{Wei}). Moreover, in the case of a cost function $J$ of the form (\ref{P2 cost obj}), one can  prove   the uniqueness of the optimal bilinear control under some constraint relaying $T$ and $y_0$ \cite{Bradly,Boukhari R,Boukhari 2017}.
\end{rem}

\section{Conclusion}

In this work, we studied the question of quadratic optimal control with endpoint constraint for bilinear systems. The optimal control is characterized via a set of unconstrained minimization problems,  then it is expressed as a time varying feedback for  commutative bilinear systems. The obtained results are applied to parabolic and hyperbolic PDE. As an interesting continuation of the present work, one can consider  the same questions for unbounded control operators,  such as the case of Fokker Planck equation \cite{fokker}.\\

{\bf Conflict of interest statement.}\\
On behalf of all authors, the corresponding author states that there is no conflict of interest.

\end{document}